\magnification=\magstep1
\hsize=15truecm
\vsize=23.5truecm

\overfullrule=0pt
\input mssymb
\def\sqr#1#2{{\vcenter{\hrule height.#2pt
      \hbox{\vrule width.#2pt height#1pt \kern#1pt
         \vrule width.#2pt}
       \hrule height.#2pt}}}
\def\op{\mathchoice\sqr34\sqr34\sqr{2.1}3\sqr{1.5}3}

\def\Jbeop{$\op $}

\def\eop{\ifmmode\op\else\Jbeop\fi }
\newdimen\refskip
\refskip=40pt                            
\def\ref #1 #2\par{\par\noindent\rlap{#1}\hskip\refskip
\hangindent\refskip #2}
\def\kohta #1 #2\par{\par\indent\rlap{\rm(#1)}\hskip\refskip
\hangindent40pt #2\par}
\vskip2truecm
\centerline{{\bf THE SUBSPACE PROBLEM}}
\centerline{{\bf FOR WEIGHTED INDUCTIVE LIMITS}}
\centerline{{\bf OF SPACES OF HOLOMORPHIC FUNCTIONS}}
\bigskip
\vskip1truecm
\centerline{{\bf Jos\'{e} Bonet and Jari Taskinen}}
\footnote{}{Mathematics subject classification: primary 46E10, secondary
46A13, 32A37}
\vskip1truecm
\bigskip
The aim of the present article is to solve in the negative a well-known
open 
problem raised by Bierstedt, Meise and Summers in [BMS1] (see also [BM1]).
We construct a countable 
inductive limit of weighted Banach spaces of holomorphic functions, 
which is not a topological subspace of the 
corresponding weighted inductive limit of spaces of continuous functions.
As a 
consequence the topology of the weighted inductive limit of spaces of 
holomorphic functions cannot be described by the weighted sup-seminorms
given 
by the maximal system of weights associated with the sequence of 
weights defining the inductive limit.
The main step of our construction shows that a certain sequence 
space is isomorphic to a complemented subspace of a weighted space of 
holomorphic functions. To do this we make use of a special sequence of outer 
holomorphic functions and of the existence of radial limits of holomorphic 
bounded functions in the disc.
\bigskip
Weighted spaces and weighted inductive limits of spaces of holomorphic 
functions on open subsets of ${\Bbb C}^N$ $(N \in {\Bbb N})$ arise in many 
fields like linear partial differential operators, convolution equations, 
complex and Fourier Analysis and distribution theory. Since the structure of 
general locally convex inductive limits is rather complicated and many 
pathologies can occur, the applications of weighted inductive limits have
been 
restricted. The reason was that it did not seem possible to 
describe the inductive limit, 
its topology, and in particular a fundamental system of seminorms 
in a way that 
permits direct estimates and computations. In the theory of 
Ehrenpreis [E] of "analytically uniform spaces", he needed that the topology 
of certain weighted inductive limits of spaces of entire functions, which are 
the Fourier-Laplace transforms of spaces of test functions or 
ultradistributions, has a fundamental system of weighted sup-seminorms. 
Berenstein and Dostal [BD] reformulated the problem in a more general setting 
and used the term "complex representation". This corresponds exactly with the 
term "projective description" used by Bierstedt, Meise and Summers [BMS1] 
which is the one we will also utilize in this paper. In [BMS1] it was proved 
that countable weighted inductive limits of Banach spaces of holomorphic 
functions on arbitrary open subsets $G$ of ${\Bbb C}^N$ admit such a canonical 
projective description by weighted sup-seminorms whenever the linking maps 
between the generating Banach spaces are compact. This theorem extended 
previous work by B. A. Taylor [T] with a more functional analytic approach and
was very satisfactory from the point of view of applications. It remained open 
whether the projective description theorem continued to hold for weighted 
inductive limits of spaces of holomorphic functions without any restriction on 
the linking maps. This problem is solved here.
\bigskip

{\bf 1. Notation and preliminaries.} 
\bigskip
All the vector spaces are defined over the complex scalar field ${\Bbb C}$. 
We denote by ${\Bbb R}^+$ (resp. ${\Bbb R}^+_0$)
the space of strictly positive reals (resp. ${\Bbb R}^+ \cup \{ 0 \}$).

Let $V = (v_k )_{k=1}^{\infty}$ be a
decreasing sequence of continuous strictly positive weight functions defined
on 
an open subset $G$ of
${\Bbb C}^N$, $N \in {\Bbb N}$. 
We denote 
by ${\cal V} C (G) $ and
${\cal V} H (G)$ 
the inductive limits ${\rm ind}_k C v_k (G)$ and ${\rm ind}_k H v_k (G)$,
where $C v_k (G)$ (respectively, $H v_k (G)$)  denotes the Banach space
$$\{ f: G \to {\Bbb C} \ {\rm continuous \ (resp.
holomorphic) } \ |$$
$$  p_{v_k}(f) := \sup_{z \in G} v_k(z) | f(z) |
< \infty          \}.
$$
The canonical embedding ${\cal V} H (G) \hookrightarrow {\cal V} C (G)$
is continuous,
and it is a well known open problem, if the topologies of ${\cal V} H
(G)$
and ${\cal V} C (G)$ coincide on ${\cal V} H (G)$. (See [BM1], Section
1 or [BiBo4], Section 4, Problem 5.) This is a particular case of the 
so-called subspace problem for locally convex inductive limits. 

In order to describe the topology of the weighted inductive limits 
${\cal V} C (G) $ and ${\cal V} H (G) $, Bierstedt, Meise and Summers [BMS1] 
introduced the system of weights 
$\overline V$, associated with the sequence $V$,
$$
\overline V = \{ \overline v : G \to  {\Bbb R}^{+} \ {\rm continuous}
\ | \ \forall k \in
{\Bbb N} \ \exists C_k > 0 \ {\rm such \ that} \ \overline v \le C_k v_k
\}.
$$
The projective hulls $C \overline V (G)$ (resp. $H \overline V (G)$) 
of ${\cal V} C (G) $ (resp. ${\cal V} H (G) $)  is the 
locally convex space
$$\{ f: G \to {\Bbb C} \ {\rm continuous \ (resp.
holomorphic) } \ | $$
$$ p_{\overline v}(f) := \sup_{z \in G}   \overline v (z)
| f(z) | < \infty   \ {\rm for \ all }\ \overline v \in \overline V     \}.
$$
endowed with the locally convex topology defined by the seminorms 
$p_{\overline v}$ as ${\overline v}$ varies in ${\overline V}$.
Clearly the inclusions ${\cal V} C (G) \hookrightarrow  C \overline V (G)$
and
${\cal V} H (G) \hookrightarrow  H \overline V (G)$ are continuous. In [BMS1] 
it was proved that ${\cal V} C (G) = C \overline V (G)$ and
${\cal V} H (G) = H \overline V (G)$ hold algebraically and that the two 
spaces in each equality have the same bounded sets. Moreover one of the main 
results in [BMS1] shows that if $V$ satisfies condition $(S)$
\bigskip
(S) \ for all $k$ there is $l$ such that $v_l/v_k$ vanishes at infinity on $G$
\bigskip
\noindent
then ${\cal V} H (G) = H \overline V (G)$ holds topologically and 
${\cal V} H (G)$ is a topological subspace of ${\cal V} C (G)$. In [BM2], [Ba] 
and [BiBo3] the topological identity ${\cal V} C (G) = C \overline V (G)$ was 
characterized in terms of a condition $(D)$ on the sequence $V$. We present
here an
example showing that if condition $(S)$ does not hold the space 
${\cal V} H (G)$ need not be a topological subspace of 
${\cal V} C (G)$ and ${\cal V} H (G) = H \overline V (G)$ need not hold 
topologically.

In the construction of our example we need weighted inductive limits of spaces 
of sequences on ${\Bbb N}$. We recall the notations from [BMS2]. We will 
denote here by $\Lambda = (\lambda_k )_{k=1}^{\infty}$ a decreasing sequence 
of strictly positive weights on ${\Bbb N}$, ${\lambda_{nk}} := \lambda_k (n)$ 
for $k,n \in {\Bbb N}$. The corresponding weighted inductive limit is denoted 
by $k_{\infty} = {\rm ind}_k l_{\infty}(\lambda_k)$. The system of 
weights associated with $ \Lambda $ is denoted by $\overline{\Lambda}$ and 
$\overline{\lambda} \in \overline{\Lambda}$ if 
and only if $\overline{\lambda}(n) > 0$ 
for every $n \in {\Bbb N}$ and for every $k \in {\Bbb N}$ there is $C_k > 0$
with $\overline{\lambda} \le 
C_k \lambda_k$ on ${\Bbb N}$. The projective hull of the inductive limit 
$ k_{\infty}$ is denoted by $K_{\infty} $, and it is the space 
$$\{ x =(x_n) \ |  p_{\overline {\lambda}}(x) := \sup_{n \in {\Bbb N}}   
\overline \lambda (n)
| x_n | < \infty   \ {\rm for \ all }\ \overline {\lambda} \in \overline
{\Lambda}     \}.
$$
The spaces $K_{\infty} $ and $k_{\infty} $ always coincide 
algebraically and they have the same bounded sets, but there are examples of 
sequences $\Lambda$ such that $K_{\infty} $ and $k_{\infty} $ do not 
coincide topologically, $K_{\infty} $ has bounded sets which are not 
metrizable and it is not bornological. See [BMS2], [BiBo1], [K] and [V]. We 
refer to [BiBo4] for a survey article on spaces of type ${\cal V} C (X)$.
\bigskip

{\bf 2. Main construction.} 
\bigskip
In this section we construct a sequence of weights $W=(w_k)^{\infty}_{k=1}$
on 
an open bounded set $G_1$ of ${\Bbb C}$ such that the projective hull 
$H \overline W (G_1)$ contains a complemented subspace isomorphic to a space 
of sequences $K_{\infty}$ which is not bornological. Consequently, the space 
$H \overline W (G_1)$ is not bornological and, hence, it does not coincide 
topologically with ${\cal W} H (G_1)$.

We first select a decreasing sequence $\Lambda = (\lambda_k )_{k=1}^{\infty}$ 
of strictly positive functions $\lambda_k(n)=\lambda_{nk}$, 
$n,k \in {\Bbb N}$, on ${\Bbb N}$ such that 
$1/n^2<\lambda_{nk}\le 1$ for all $n$ and $k$, 
the corresponding space $K_{\infty}$ is not 
bornological and it contains bounded sets which are not metrizable. 
For example combine [K], Section 31.7 with [BM1], Theorem 9 and [BiBo1].
In this case the system of weights $\overline {\Lambda}$ associated with 
$\Lambda$ satisfies that for each $\overline {\mu} \in \overline {\Lambda}$ 
there are $\overline {\lambda} \in \overline {\Lambda}$ and $C>0$ with 
$\overline {\mu} \le C \overline {\lambda}$ and
$1/n^2<\overline{\lambda}(n) \le 1$ for all $n \in {\Bbb N}$. Indeed, 
given $\overline {\mu} \in \overline {\Lambda}$ , we select $c_k>0$ such that 
$\overline {\mu} \le \inf_k c_k\lambda_k$. We put $d_k = \max(c_k,1)$ for all 
$k$ and we set $\overline {\lambda} = \min(\inf_k d_k\lambda_k,\lambda_1) \in 
\overline {\Lambda}$. Accordingly $K_\infty$ has a fundamental system of 
seminorms P given by multiples of elements 
$\overline {\lambda} \in \overline {\Lambda}$ satisfying
$1/n^2<\overline{\lambda}(n) \le 1$ for all $n \in {\Bbb N}$.

We denote $G_1=\{ z\in\Bbb C\mid 1/2<|z|<1,\
0<\arg z<\pi\}$ and
we define the system
$W=(w_k)^{\infty}_{k=1}$ of weight functions on $G_1$ by
$$w_k(re^{i\theta})=\hat w_k(\theta ),$$
where $\hat w_k:[0,\pi [ \to \Bbb R^+$ satisfies
$$\hat w_k(\theta )=\lambda_{nk}$$
for $\theta\in I_n:=[\theta_n-1/(2^5n^2),\ \theta_n+1/(2^5n^2)],\
\theta_n:=1/(2n)$, for all $n\in\Bbb N$,
$$\hat w_k(\theta_n+1/(2^4n^2))=\hat w_k(0)=\hat w_k(\pi)=1,$$
for all $n$, and $\hat w_k$ is extended affinely for other $\theta$.
\medskip

Now we define a sequence of elements of $H(G_1)$ which will be essential in 
our construction. For all $n\in\Bbb N$, 
let $\varepsilon_n>0$ be a number satisfying
$$\varepsilon_n<2^{-n-16}n^{-6};$$
because of the choice of $(\lambda_{nk})$ we have, in particular,
$$\varepsilon_n<2^{-n-16}n^{-4}\lambda_{nk}.$$
For all $n$, let $e_n$ be an analytic function on the disc defined by
$$e_n(z)=\exp \big(
{1\over2\pi}\int\limits_0^{2\pi}{e^{i\theta}+z\over
e^{i\theta}-z}\log\varphi_n(\theta )d\theta ),$$
where $\varphi_n:[0,2\pi]\to\Bbb R^+$ is the measurable function
$$\eqalign{\varphi_n(\theta ):=\cases{1\quad {\rm for}\ \theta\in
J_n:=[\theta_n-\varepsilon_n,\
\theta_n+\varepsilon_n]\cr\noalign{\vskip4pt}\varepsilon_n2^{-m-4}\quad
{\rm for}\ \theta\in J_m,\ m\ne
n\cr\noalign{\vskip4pt}\varepsilon_n\quad {\rm for\ other}\ \theta
.\cr}\cr}$$
In fact, by $[R]$, Theorem 17.16, $e_n \in H^{\infty}$ and 
$$\mid e^*_n(e^{i\theta})\mid =\varphi_n(\theta ),$$
where $e^*_n(e^{i\theta}):=\lim\limits_{r\to 1}e_n(re^{i\theta})$, holds
for a.e. $\theta\in [0,2\pi ]$. We also denote by $e_n$ the restrictions of 
$e_n$ to $G_1$.

For $n\in\Bbb N$ we denote
$$D_n:=\{ z\in G_1\mid |z-e^{i \theta_n} |<1/(50n^2)\}, \ \
C_n:=G_1 \setminus D_n.$$
We have $D_n \subset \{ z \in G_1\mid z=re^{i \theta} , \ \theta \in I_n \}$ 
and, moreover, 
$|e^{i\theta}-z|>1/(2^6n^2)$ for $\theta\in J_n$ and $z\in C_n$.
Since
$$\int\limits_0^{2\pi}(1-|z|^2)/|e^{i\theta}-z|^2d\theta =2\pi ,$$
we can apply
the Jensen inequality ($[R]$, Theorem 3.3; take $\exp$ for
the convex function
and $(1-|z|^2)/(2\pi |e^{i\theta}-z|^2)d\theta$ for the probablity
measure) to get for $z\in C_n$
$$
\eqalignno{&|e_n(z)|\le\exp \big( {1\over 2\pi}{\rm
Re}\int\limits_0^{2\pi}{e^{i\theta}+z\over
e^{i\theta}-z}\log\varphi_n(\theta )d\theta )\cr\noalign{\vskip4pt}
&= \exp\big( {1\over 2\pi}\int\limits_0^{2\pi}{1-|z|^2\over
|e^{i\theta}-z|^2}\log\varphi_n(\theta )d\theta )\cr\noalign{\vskip4pt}
&\le {1\over 2\pi}\int\limits_0^{2\pi}{1-|z|^2\over
|e^{i\theta}-z|^2}\varphi_n(\theta )d\theta\cr\noalign{\vskip4pt}
&\le {1\over 2\pi}\int\limits_{J_n}{1-|z|^2\over
|e^{i\theta}-z|^2}d\theta +{\varepsilon_n\over 2\pi}\int\limits_{[0,2\pi
] \setminus J_n}{1-|z|^2\over |e^{i\theta}-z|^2}d\theta }$$
$$ \le \pi^{-1}\varepsilon_n2^{12}n^4+\varepsilon_n\le
2^{-4-n}.$$

Analogously, one can show $ |e_n(z)|\le 1$ for all $z\in G_1$. 
\bigskip
In the proof of our main result in this section we need the following 
technical lemma. It shows that 
given an arbitrary weight function $\overline w' \in \overline W$ we
can choose a dominating weight function on $G_1$ which has
some specific properties. 
\bigskip
{\bf Lemma 1.} {\sl Given a weight function $\overline w' \in
\overline W$ we can find a weight function $\overline w \in
\overline W$ with the following properties:

$1^{\circ}.$ There exists $C > 0$ such that $C \overline w'
\le \overline w \le 1$.

$2^{\circ}.$ If $\overline w'' \in \overline W$ is defined as
$w_k $  except that  $\lambda_{nk}$ is replaced by
$1/n^2$,  we have $\overline w'' \le \overline w$.

$3^{\circ}.$ The weight  $\overline w$ is constant on every
$D_n$  so
that
$$
\overline \lambda (n) : {\Bbb N} \to {\Bbb R}^{+}, \ \ \overline
\lambda (n) := \overline w (z) , \ z \in D_n,
$$
satisfies $1/ n^2 \le \overline \lambda (n) $ 
for all $k $ and $n$ and $\overline \lambda \in \overline \Lambda$.}
\bigskip
Proof. Let $\varrho (n) := \max\{ 1/n^2 , \sup \{ \overline w'(z)
|\  {\rm arg}(z) \in I_n \} \}$ and define $\overline w^{(1)}$ as
$w_k$  but replace $\lambda_{nk} $ by $\varrho (n)$.
Now it is easy to see that the weight $\overline
w$, defined by
$$\overline w(z) = \min \{ w_1 (z) , \max \{ \overline w^{(1)} ,
\overline w'(z) \} \}
$$
for $z \in G_1$, has all the desired properties; the property
$3^{\circ}$ follows from the facts that $\overline
w^{(1)}$ and $w_1$ are constants on $D_n$ and $\overline w^{(1)}
\ge \overline w' $ on $D_n$. \quad \eop
\bigskip
Our next lemma is essentially known.
\bigskip
{\bf Lemma 2.} {\sl Let $E$ and $F$ be complete locally convex spaces. Let 
$\psi : E \rightarrow F$ and $\phi : F \rightarrow E$ be continuous linear 
maps such that $\phi \psi : E \rightarrow E$ satisfies the following 
condition: there is a fundamental system of seminorms $P$ on $E$ and there is 
$0 < \delta < 1$ such that 
$$p((\phi \psi - id_E)x) \le \delta p(x) \ \ \ \forall x \in E \ \forall p 
\in P.$$
Then $E$ is isomorphic to a complemented subspace of $F$.}
\bigskip
Proof. We put $B:= \phi \psi - id_E$ and we define $A: E \rightarrow E$ by 
$$Ax := \sum_{n=0}^{\infty} (-1)^n B^n x \ , \  x \in E.$$
Then $A$ is a well defined continuous linear operator on $E$. Indeed, for $x 
\in E$, the series $\sum_{n=0}^{\infty} (-1)^n B^n x$ is absolutely summable
 in 
$E$ and, for $p \in P, x \in E$, we have
$$p(\sum_{n=0}^{\infty} (-1)^n B^n x) \le \sum_{n=0}^{\infty} p(B^n x) \le 
\sum_{n=0}^{\infty} \delta^n p(x) \le C p(x).$$
Moreover $\phi (\psi A) = id_E$. Indeed
$$\phi (\psi A) = (\phi \psi) A = (id_E + B)\sum_{n=0}^{\infty} (-1)^n B^n = 
id_E.$$
This implies that $(\psi A) \phi$ is a projection on $F$ whose image is 
isomorphic to $E$ (see e.g.\ [H], pp.\ 122-123). \quad \eop
\bigskip
{\bf Theorem 3.} {\sl The space $H \overline W (G_1)$ contains a complemented
subspace isomorphic to the non-bornological space $K_{\infty}$.In particular 
$H \overline W (G_1)$ does not coincide topologically with the weighted 
inductive limit ${\cal W} H(G_1)$.}
\bigskip
Proof. We construct continuous linear maps $\psi : K_{\infty} \rightarrow 
H \overline W (G_1)$ and $\phi : H \overline W (G_1) \rightarrow K_{\infty}$ 
satisfying the assumptions of lemma 2.

First define $\psi : K_{\infty} \rightarrow H \overline W (G_1)$ by 
$\psi(a):= \sum_{n=1}^{\infty} a_n e_n $ for 
$a = (a_n)_{n=1}^{\infty}   \in K_{\infty}$. To see that $\psi$ is well
defined
and continuous, we fix $\overline w' \in \overline W$, and we select 
$\overline w \in \overline W$ and $\overline \lambda \in \overline \Lambda$
as 
in Lemma 1. 
If $(a_n)^{\infty}_{n=1}$ is a sequence of
scalars such that $\sup\limits_n\overline{\lambda}(n)|a_n|=1$, we have
$| a_n | \le n^2$ for all $n$. Every $z_0 \in G_1$ has a neighbourhood
$U$
which intersects at most one of the sets $D_n$. It follows
from the estimates of $| e_n(z) |$ established after the definition of 
$e_n$ that $\sum a_n e_n (z) $
converges uniformly for $z \in U$ and
thus defines a holomorphic function of $G_1$.
Moreover,
denoting
$D:=\bigcap\limits^{\infty}_{n=1}C_n$, by the choice of $(a_n)$, we have 
$$\eqalignno{&p_{\overline w}(\sum\limits^{\infty}_{n=1}a_ne_n)
\cr\noalign{\vskip4pt}
& \le \sup\limits_{z\in
G_1}\overline w (z)\sum\limits^{\infty}_{n=1}
\overline{\lambda}(n)^{-1}|e_n(z)|
\cr\noalign{\vskip4pt}
&=\sup\limits_{m\in \Bbb N}\sup\limits_{z\in
D_m}(\overline w  (z)\sum_{\scriptstyle n\in\Bbb N,\atop\scriptstyle n\ne
m}\overline \lambda (n)^{-1} |e_n(z)|
+\overline w (z) \overline \lambda (m)^{-1} |e_m(z)|)
\cr\noalign{\vskip4pt}
&+\sup\limits_{z\in
D}\overline w (z)\sum\limits^{\infty}_{n=1}
\overline \lambda(n)^{-1} |e_n(z)| \le }
$$
$$\sup\limits_{m\in\Bbb
N}(\sum\limits^{\infty}_{n=1}\overline  \lambda(n)^{-1}
2^{-4-n} +
\overline \lambda({m}) \overline \lambda ({m})^{-1})
$$
$$+\sum\limits_{n=1}^{\infty}\overline \lambda ({n})^{-1}2^{-4-n}
\le
3 = 3\sup\limits_n\overline{\lambda}(n)|a_n|.$$
This shows that the map $\psi$ is continuous.
\bigskip
To define $\phi : H \overline W (G_1) \rightarrow K_{\infty}$ we need radial 
values of elements of $H \overline W (G_1)$. We fix $f \in 
H \overline W (G_1)$. There is $k \in {\Bbb N}$ such that $f \in H w_k
(G_1)$. 
Given $n \in {\Bbb N}$, since the weight $w_k$ is constant in 
$ \{ r e^{i \theta} : \theta \in I_n , 1/2 < r < 1 \} $, it follows that 
$f^*(e^{i\theta})
=\lim\limits_{r\to 1}f(re^{i\theta})$ exists a.e. for $\theta \in J_n$. See 
[HF], pp.\ 34 and ff.\ or [D], pp.\ 170. Accordingly we define 
$f^*(e^{i\theta})$ a.e. $\theta \in J_n$, which is an element of 
$L^{\infty}(J_n)$. Observe that the radial limits of $e_n$ in each $J_n$ are 
the restriction of the ones of $e_n$ in the disc, and that, if $(a_n) \in 
K_{\infty}$, it follows from the inequalities established in the first part
of 
this proof, that $(\sum a_n e_n)^*(e^{i\theta}) = \sum a_n 
e_n^*(e^{i\theta})$. We set 
$\chi_n(\theta ):=e^{-i\arg e^*_n(e^{i\theta }) }$ and we define 
$\phi : H \overline W (G_1) \rightarrow K_{\infty}$ by 
$$
\phi(f):=
((2 \varepsilon_n)^{-1} \int\limits_{J_n}f^*(e^{i\theta})
\chi_n(\theta )d\theta )_{n \in {\Bbb N}}
$$ 
We first check that $\phi$ is well defined and continuous. 
Given $\overline{\lambda} \in \overline{\Lambda}$ 
with $1/n^2\le
\overline{\lambda}(n)\le 1$ for all $n$, we
define the weight $\overline w \in \overline W$
as $w_k$, but replacing
$\lambda_{nk}$ by $\overline \lambda (n)$.
If $f\in H\overline W (G_1)$ and $n \in {\Bbb N}$ we have
$$
\overline \lambda(n) (2 \varepsilon_n)^{-1}| \int\limits_{J_n} 
f^*(e^{i\theta}) \chi_n(\theta) d\theta | \le
\overline \lambda (n) \sup\limits_{J_n} | f^*(e^{i\theta}) | \le 
\sup_{z \in G_1} \overline w(z) |f(z)|.
$$
This shows  $\phi(f) \in K_{\infty}$ and the continuity of $\phi$.
\bigskip
It remains to show that $\phi \psi - id_{K_{\infty}}$ satisfies the condition 
in lemma 2. First observe that for all $n \in {\Bbb N}$
$$
(2 \varepsilon_n)^{-1} \int\limits_{J_n}e_n^*(e^{i\theta})
\chi_n(\theta )d\theta = 1.
$$
On the other hand, for $n \in {\Bbb N}$ fixed, $| e_j^*(e^{i\theta})| =
\varepsilon_j 2^{-n-4}$ for all $j \in {\Bbb N} , j \ne n$ , a.e. $\theta \in 
J_n$ and the series $\sum a_j \varepsilon_j$ converges absolutely for $(a_n) 
\in K_{\infty}$. Indeed, select 
$\overline{\lambda} \in \overline{\Lambda}$ 
with $1/n^2\le
\overline{\lambda}(n)\le C$ for all $n \in {\Bbb N}$. We have $S:= \sup_j 
\overline \lambda(j) |a_j| < \infty$ and $\varepsilon_j < 2^{-j-16} j^{-6} < 
2^{-j-16} j^{-4} \overline \lambda(j)$, then $\sum |a_j| \varepsilon_j < S 
\sum 2^{-j-16}.$ In particular
$$ \sum |a_j| \varepsilon_j \le (1/8) \sup_n \overline \lambda(n) |a_n| \ \ \ 
\forall \overline{\lambda} \in \overline{\Lambda}.$$
Moreover
$$(\phi \psi - id_{K_{\infty}})((a_n)) = ((2 \varepsilon_n)^{-1}
\sum\limits_{j \ne n} a_j 
\int\limits_{J_n} 
e_j^*(e^{i\theta})
\chi_n(\theta )d\theta)_{n \in {\Bbb N}}.$$ 
If $\overline{\lambda} \in \overline{\Lambda}$ satisfies 
$1/n^2\le
\overline{\lambda}(n)\le 1$ for all $n \in {\Bbb N}$ and $a \in K_{\infty}$ is 
such that $p_{\overline \lambda}(a) = \sup_n \overline \lambda (n) |a_n| = 1$ 
we have, for $n \in {\Bbb N}$,
$$
| \overline \lambda(n) (2 \varepsilon_n)^{-1} \sum\limits_{j \ne n} a_j 
\int\limits_{J_n} 
e_j^*(e^{i\theta})
\chi_n(\theta )d\theta | \le 
\overline \lambda(n) (2 \varepsilon_n)^{-1} \sum\limits_{j \ne n} | a_j |
\int\limits_{J_n} 
| e_j^*(e^{i\theta})
\chi_n(\theta ) | d\theta \le
$$
$$
\overline \lambda(n) \sum\limits_{j \ne n} | a_j | 
\sup\limits_{\theta \in J_n} 
| e_j^*(e^{i\theta}) | \le 2^{-n-4} 
\overline \lambda(n) \sum\limits_{j \ne n} | a_j | \varepsilon_j \le 
1/128.
$$
Therefore
$$
p_{\overline \lambda}((\phi \psi - id_{K_{\infty}})(a)) \le 
(1/128) p_{\overline \lambda}(a).
$$
Since the multiples of the weights
$\overline{\lambda} \in \overline{\Lambda}$ with 
$1/n^2\le
\overline{\lambda}(n)\le 1$ for all $n \in {\Bbb N}$ 
form a fundamental system of seminorms of the space $K_{\infty}$, 
the conclusion follows from lemma 2. \quad \eop

\bigskip

{\bf 3. The subspace problem.}
\bigskip
In this section we denote 
$G=G_1\times \Bbb C\subset \Bbb C^2$ and we construct a decreasing 
sequence  $V=(v_k)^{\infty}_{k=1}$ of weight functions on $G$ such that 
${\cal V} C (G) = C \overline V (G)$ holds topologically, but 
${\cal V} H (G)$ is not a topological subspace of ${\cal V} C (G)$. 
Moreover the projective hull $H \overline V (G)$ is not even a (DF)-space.

If $z_1 \in G_1$, we write $d(z_1)$ to denote the distance of $z_1$ to the 
complement of $G_1$.

We define the
system $V=(v_k)^{\infty}_{k=1}$ of weight functions on $G$ by
$$v_k(z_1,z_2)=w_k(z_1)u_k(z_1,|z_2|),$$
where $u_k:G_1\times\Bbb R_0^+\to\Bbb R^+$ is defined by
$$\eqalign{u_k(z_1,t)=\cases{(1 +t)^{-{k-1\over 2k}},\quad t\ge
k+1\cr\noalign{\vskip4pt} (1 + {1\over d(z_1)}+t)^{-{k-1\over 2k}},\quad
t\le k\cr}\cr}$$
and, for each fixed $z_1,\ u_k(z_1,t)$ is extended affinely for
$k<t<k+1$. It is easy to see that the functions $v_k$ are continuous on
$G$.

Bierstedt and Meise [BM2] introduced the following {\it condition} $(M)$ on 
the sequence $V=(v_k)^{\infty}_{k=1}$ : for each $k \in {\Bbb N}$ and each 
subset $Y$ of $G$ which is not relatively compact, there exists $k' = k'(k,Y) 
>k$ with $\inf\limits_{y \in Y} v_{k'}(y)/v_k(y) = 0$. They proved that this 
condition is equivalent to the fact that $C \overline V (G)$ induces the 
compact open topology on each bounded subset and that condition 
$(M)$ implies the 
topological identity 
${\cal V} C (G) = C \overline V (G)$. Moreover, if $V$ satisfies $(M)$, then 
$H \overline V (G)$ is a Montel space. It was an open problem (see also [Bi]) 
whether ${\cal V} H (G) $ is a Montel space when $V$ satisfies $(M)$. In fact, 
if ${\cal V} H (G) $ is a Montel space, then 
${\cal V} H (G) = H \overline V (G)$ holds topologically 
by a direct application of the Baernstein
open mapping lemma as in [BMS1].
\bigskip
{\bf Proposition 4.} {\sl The sequence $V=(v_k)^{\infty}_{k=1}$ satisfies
condition 
$(M)$. Consequently, ${\cal V} C (G) = C \overline V (G)$ holds algebraically 
and topologically, and $H \overline V (G)$ is a Montel space with metrizable 
bounded sets.}
\bigskip
Proof. Let $k \in {\Bbb N}$ be given and let $Y$ be a subset of $G$ which is 
not relatively compact. We have two possibilities: either
$$\eqalign{&{\rm (i)}\ \exists 
(z^{(m)})=((z_1^{(m)},z_2^{(m)})) \subset Y  \ : \ 
\sup\limits_m|z_2^{(m)}|=\infty,\ {\rm
or}\cr\noalign{\vskip4pt}
&{\rm (ii)}\ \exists 
(z^{(m)})=((z_1^{(m)},z_2^{(m)})) \subset Y : | z_2^{(m)} | \le M  \ 
{\rm and} \  
\inf\limits_m\ d(z_1^{(m)})=0.\cr}$$
In the first case it follows easily from the definition of $u_k$ that,
taking $k' = k+1$,
$$\eqalign{&\sup\limits_{m\in\Bbb N}{v_k(z^{(m)})\over
v_{k+1}(z^{(m)})}\ge \sup\limits_{m\in\Bbb
N}{u_k(z_1^{(m)},|z_2^{(m)}|)\over
u_{k+1}(z_1^{(m)},|z_2^{(m)}|)}\cr\noalign{\vskip4pt}
&\ge\sup_{\scriptstyle m\in\Bbb N,\atop\scriptstyle
|z_2^{(m)}|>k+2}|z_2^{(m)}|^{-{k-1\over 2k}+{k\over 2(k+1)}}=\infty .\cr}$$
In case (ii), we choose $k'>M$ and $k' > k+1$ to get
$$\sup\limits_{m\in\Bbb N}{v_k(z^{(m)})\over v_{k'}(z^{(m)})}\ge
\sup\limits_{m\in\Bbb N}\big( {1\over d(z_1^{(m)})}+1+M)^{-{k-1\over
2k}}\cdot \big( {1\over d(z_1^{(m)})}\big)^{{k\over 2(k+1)}}=\infty
.\quad\eop$$
\bigskip
It is very easy to see that every $f\in H\overline V(G)$ 
is constant with respect to
the second variable. Indeed, if $f\in H\overline V(G)$, 
there are $k\in\Bbb N$ and $C>0$ such that 
$p_{v_k}(f)\le C,$
so that for every fixed $z_1\in G_1$
$$\eqalignno{&C\ge\sup\limits_{z_2\in\Bbb
C}w_k(z_1)u_k(z_1,|z_2|) | f(z_1,z_2) | \cr\noalign{\vskip4pt}
&\ge w_k(z_1)\sup\limits_{z_2\in\Bbb C}\big(1 + {1\over 
d(z_1)}+|z_2| \big)^{-1/2}| f(z_1,z_2) | .\cr}$$
Now it is an elementary fact of complex analysis that a holomorphic
$g:\Bbb C\to \Bbb C$ satisfying $\sup\{(|z|+c_0)^{-1/2} | g(z) |  \
\mid z\in\Bbb
C\} <\infty$ for some constant $c_0$, must be constant. Accordingly 
$f$ must be constant  with respect to $z_2$.
Now we define $A: H\overline V (G) \rightarrow H\overline W (G_1)$ by
$Af(z_1) = f(z_1,0)$. To show that $A$ is well defined, we observe that 
$$
\eqalignno{&p_{w_k}(Af):=\sup\limits_{z\in G_1}w_k(z)|f(z,0)|\le
\sup\limits_{z\in
G_1}w_k(z)C_ku_k(z , k+1)| f(z,0) | \cr\noalign{\vskip4pt}
&\le C_kp_{v_k}(f)
&\cr}
$$
for all $f\in Hv_k(G)$
and $C_k:=(k+2 )^{(k-1)/(2k)}$.
\bigskip
Given $g \in H\overline W (G_1)$, we define $\overline g : G \rightarrow
{\Bbb
C} $ by $\overline g(z_1,z_2) = g(z_1)$ for all $(z_1,z_2) \in G$. To show 
that $\overline g \in 
H\overline V (G)$, we fix $k \in {\Bbb N}$ with $g \in H w_k (G_1)$, then we 
have the estimate 
$$
\sup\limits_{(z_1,z_2) \in G} v_k(z_1,z_2) | \overline g(z_1,z_2) | \le 
\sup\limits_{z_1 \in G_1} w_k(z_1) | g(z_1) | ,
$$
\noindent
since $0 \le u_k \le 1$. This shows that A is bijective and that 
$A^{-1} : H w_k (G_1) \rightarrow H v_k (G)$ is continuous for every $k \in 
{\Bbb N}$. By the closed graph theorem for (LB)-spaces, this also yields that 
$A : {\cal V} H (G) \rightarrow {\cal W} H (G_1)$ is a topological 
isomorphism. Moreover, $A^{-1} : H \overline W (G_1) \rightarrow 
H \overline V (G)$ is continuous. Indeed, if $\overline v \in \overline V$ is 
given, we define 
$\overline w (z_1) = \sup\limits_{z_2 \in {\Bbb C}} \overline v
(z_1 , z_2) $ for $z_1 \in G_1$. We  have for all $k$
$$
\overline w (z_1) = \sup\limits_{z_2 \in {\Bbb C}} \overline v
(z_1 , z_2) \le \sup\limits_{z_2 \in {\Bbb C}} C_k v_k ( z_1 , z_2)
\le C_k w_k(z_1),
$$
hence $\overline w \in \overline W$. Moreover,
$$
p_{\overline v } (f) = \sup\limits_{(z_1,z_2 ) \in G}
\overline v (z_1 , z_2) | f (z_1,0 ) | $$
$$\le \sup\limits_{z_1 \in G_1}
 \{ | f (z_1,0 ) | \sup\limits_{z_2  \in {\Bbb C} }
\{ \overline v  ( z_1 , z_2) \} \}
= p_{ \overline w } (Af) .
$$
On the other hand $A : H \overline V (G) \rightarrow H \overline W (G_1)$ is 
not continuous. In fact, $H \overline W (G_1)$ and $H \overline V (G)$ can not 
be isomorphic, since the first one contains a complemented subspace isomorphic 
to $K_{\infty}$ by theorem 3; 
hence it contains bounded sets which are not metrizable; 
while every bounded subset of $H \overline V (G)$ is metrizable by proposition 
4.
\bigskip
{\bf Theorem 5.} {\sl The space $H \overline V (G)$ is not bornological, 
${\cal V} H (G) = H \overline V (G)$ does not hold topologically, 
${\cal V} H (G)$ is not a topological subspace of
${\cal V} C (G)$ and it is not a Montel space. Moreover, 
$H \overline V (G)$ is not a (DF)-space.}
\bigskip
Proof. If $H \overline V (G)$ is bornological (or equivalently if 
${\cal V} H (G) = H \overline V (G)$ holds topologically), then the linear map
 $A^{-1} : H \overline W (G_1) \rightarrow {\cal V} H (G) $ is continuous. 
Consequently the identity $id = A A^{-1} : H \overline W (G_1) \rightarrow 
{\cal W} H (G_1)$ is continuous. This implies $H \overline W (G_1)$ is 
bornological; which contradicts theorem 3.

By proposition 4, ${\cal V} C (G) = C \overline V (G)$ holds topologically. 
Since $H \overline V (G)$ is clearly a topological subspace of 
$C \overline V (G)$, we conclude that 
${\cal V} H (G)$ is not a topological subspace of ${\cal V} C (G)$. If 
${\cal V} H (G)$ were Montel, we could apply directly Baernstein open mapping 
lemma (see e.g. [BPC]; 8.6.8(5)) to conclude that ${\cal V} H (G)$ would be a 
topological subspace of ${\cal V} C (G)$; a contradiction.

Finally assume $H \overline V (G)$ is a (DF)-space. Since it is a complete
(DF)-space which is Montel, $H \overline V (G)$ is bornological (cf. [BPC]; 
8.3.48). This is a contradiction. \quad \eop
\bigskip\bigskip
{\bf Acknowledgement:} The research of J.Bonet was partially supported by
DGICYT Proyecto no. PB91-0538.
Part of the research of J.Taskinen was carried
out during a visit to Paderborn  in November 1992. He wishes to thank
Prof. K.D.Bierstedt for hospitality and the Universit{\"a}t--GH Paderborn
for support. Both authors want to thank Prof. Bierstedt for many
discussions on the subject.

\bigskip\bigskip
\centerline{{\bf References:}}
\bigskip\bigskip\bigskip

\ref[Ba] Bastin, F.: On bornological spaces $C \overline V (X)$. 
Archiv Math. 53 (1989) 393-398.

\ref[BD] Berenstein, C.A., Dostal, M.A.: Analytically Uniform Spaces and their 
Applications to Convolution Equations. Springer Lecture Notes in Math. 256 
(1972).

\ref[Bi] Bierstedt, K.D.: Weighted inductive limits of spaces of 
holomorphic functions. To appear in the Proceedings of the 23rd
Annual Iranian Congress of Mathematics, April 1992..

\ref[BiBo1] Bierstedt, K.D., Bonet, J.: Stefan Heinrich's density
condition for Fr\'echet spaces and the characterization of the
distinguished K\"othe echelon spaces. Math. Nachr. 135 (1988), 149--180.

\ref[BiBo2] Bierstedt, K.D., Bonet, J.: Dual density conditions in
(DF)--spaces, I. Resultate Math. 14 (1988), 242--274.

\ref[BiBo3] Bierstedt, K.D., Bonet, J.: Dual density conditions in
(DF)--spaces, II. Bull.Soc.Roy.Sci.Li\'ege 57 (1988), 567--589.

\ref[BiBo4] Bierstedt, K.D., Bonet, J.: Some recent results on
${\cal V} C (X)$. Advances in the theory of Fr\'echet spaces, pp. 181--194,
Kluwer (1989).

\ref[BM1] Bierstedt, K.D., Meise, R.: Weighted inductive limits and
their projective descriptions. Doga Mat. 10,1 (1986), 54--82. (Special
issue: Proceedings of the Silivri Conference 1985)

\ref[BM2] Bierstedt, K.D., Meise, R.: Distinguished echelon spaces and the 
projective description of weighted inductive limits of type ${\cal V}_d C 
(X)$. Aspects of Mathematics and its Applications, pp. 169-226, Elsevier 
(1986).

\ref[BMS1] Bierstedt, K.D., Meise, R., Summers, W.H.: A projective description 
of weighted inductive limits. Trans. Amer. Math. Soc. 272 (1982) 107-160.

\ref[BMS2] Bierstedt, K.D., Meise, R., Summers, W.H.: K\"othe sets and K\"othe 
sequence spaces. Functional Analysis, Holomorphy and Approximation Theory, pp. 
27-91, North-Holland Math. Studies 71 (1982).

\ref[D] Duren P.L.: Theory of $H_p$-spaces. Academic Press (1970).

\ref[E] Ehrenpreis, L.: Fourier Analysis in Several Complex Variables. 
Interscience Tracts in Math. 17, Wiley (1970).

\ref[HF] Hoffman, K.: Banach Spaces of Analytic Functions. Prentice Hall 
(1962).

\ref[H] Horv\'ath, J.: Topological Vector Spaces and Distributions. 
Addison-Wesley (1966).

\ref[K] K\"othe, G.: Topological vector spaces, Vol.\ 1.\ Second
printing. Springer
Verlag (1983).

\ref[BPC] P\'erez Carreras, P., Bonet, J.: Barrelled Locally Convex Spaces.
North-Holland Math. Studies 131 (1987).

\ref[R] Rudin, W.: Real and complex analysis, second edition. Mc
Graw--Hill, New York (1974).

\ref[T] Taylor, B.A.: A seminorm topology for some (DF)-spaces of entire 
functions. Duke Math. J. 38 (1971) 379-385.

\ref [V] Vogt, D.: Distinguished K\"othe spaces. Math. Z. 202 (1989) 143-146.

\bigskip\bigskip\bigskip

{\bf Authors' Addresses:}
\bigskip
Departamento de Matem\'atica Aplicada, Universidad Polit\'ecnica de Valencia, 
E-46071, Valencia, SPAIN
\bigskip
Department of Mathematics, P.O Box 4 (Hallituskatu 15), SF-00014
UNIVERSITY OF HELSINKI, FINLAND

\bye